\documentclass[11pt]{article}
\usepackage{latexsym,amsmath,amsthm,amssymb,amscd,amsfonts,enumitem}
\usepackage{epsfig}
\usepackage{nicefrac}
\usepackage[all]{xy}
\usepackage{graphicx}

\usepackage[usenames]{color}

\newtheorem{question'}[varquest]{Question}
\newtheorem{thm}{Theorem}
\newtheorem{conj}{Conjecture}

\newtheorem{lemma}{Lemma}[section]
\newtheorem{remark}[lemma]{Remark}
\newtheorem{fact}[lemma]{Fact}
\newtheorem{prop}[lemma]{Proposition}

\newtheorem{theorem}[lemma]{Theorem}
\newtheorem{definition}[lemma]{Definition}
\newtheorem{claim}[lemma]{Claim}
\newtheorem{corollary}[lemma]{Corollary}

\newtheorem*{defn*}{Definition}

\newtheorem*{remark*}{Remark}

\newtheorem*{rems*}{Remarks}
\newtheorem*{thm*}{Theorem}
\newtheorem*{exs*}{Examples}

\newtheorem*{quesA'}{Question A'}
\newtheorem*{quesD'}{Question D'}

\newcommand{\R}{\mathbb{R}}
\newcommand{\C}{\mathbb{C}}
\newcommand{\D}{\mathbb{D}}

\newcommand{\N}{\mathbb{N}}
\newcommand{\Z}{\mathbb{Z}}

\newcommand{\T}{\mathbb{T}}
\renewcommand{\H}{\mathbb{H}}

\newcommand{\cc}{\mathbb{\mathcal C}}
\newcommand{\cb}{\mathbb{\mathcal B}}
\newcommand{\ca}{\mathcal A}
\newcommand{\ce}{\mathcal E}
\newcommand{\cd}{\mathcal D}

\newcommand{\cu}{\mathcal U}
\newcommand{\cm}{\mathcal M}

\newcommand{\cf}{\mathcal{F}}

\newcommand{\cj}{\mathcal{J}}

\newcommand{\id}{\textnormal{Id}\,}

\newcommand{\nbd}{neighbourhood }
\newcommand{\nbds}{neighbourhoods }
\newcommand{\fonction}[5]
{$$ 
\begin{array}{rcccl}
 #1 & : & #2 & \longrightarrow &#3 \\
    &   & #4 & \longmapsto &#5 
\end{array}
$$}

\newcommand{\dom}{\textnormal{Dom}\,}

\newcommand{\red}{\textnormal{Red}\,}

\newcommand{\rond}[1]{\overset{\vspace*{-1pt}\circ}{#1}}
\newcommand{\priv}{\backslash}
\newcommand{\lra}{\longrightarrow}
\newcommand{\hra}{\hookrightarrow}

\newcommand{\om}{\omega}
\newcommand{\eps}{\varepsilon}
\renewcommand{\phi}{\varphi}

\newcommand{\st}{\textnormal{st}}

\newcommand{\hur}{\text{Hur}}

\newcommand{\wdt}[1]{\widetilde{#1}}
\newcommand{\cqfd}{\hfill $\square$ \vspace{0.1cm}\\ }
\newcommand{\sbull}{{\tiny $\bullet$ }}

\newcommand{\im}{\textnormal{Im}\,}
\newcommand{\re}{\textnormal{Re}\,}
\renewcommand{\bar}[1]{\overline{#1}}

\newcommand{\nf}[2]{{\nicefrac{#1}{#2}}}

\newcommand{\Om}{\Omega}

\newcommand{\its}{\item[\sbull]}
\newcommand{\itss}{\vspace*{-,2cm}\item[\sbull]}

\newcommand{\psl}{{\textnormal{PSL}}}
\newcommand{\aut}{{\textnormal{Aut}}}
\newcommand{\op}{\textnormal{Op}\,}

\renewcommand{\eqref}[1]{(\ref{#1}\hspace{-,15cm})}

\newcommand{\leb}{{\cal L}\textnormal{eb}}

\def\eps{\varepsilon}

\makeatletter \@addtoreset {equation}{section}

\renewcommand\theequation
  {\ifnum \c@subsection>\z@ \arabic{section}.\arabic{subsection}.\arabic{equation}
  \else \arabic{section}.\arabic{equation} \fi}
\makeatother

\begin{document}

\title{\vspace*{-0cm} Symplectic Camel theorems and $\cc^0$-rigidity of coisotropic submanifolds}

\author{Emmanuel Opshtein}

\date{\today}
\maketitle

\begin{abstract}
This paper deals with the $\cc^0$-rigidity of the reduction of coiostropic submanifolds under the action of symplectic homeomorphism. More precisely, we exhibit several situations where a symplectic homeomorphism that takes a coisotropic submanifold to a smooth submanifold (which are then known to be coisotropic by a result of Humilière-Leclercq-Seyfaddini) abides to the non-squeezing property in the reduction. 
\end{abstract}

\section{Introduction}

This paper continues the study of the action of the symplectic homeomorphisms on smooth submanifolds as initiated in \cite{opshtein,hulese,buop}. Recall that a symplectic homeomorphism between symplectic manifolds $M$ and $M'$ is a homeomorphism that can be approximated by symplectic diffeomorphisms between $M$ and $M'$. In \cite{hulese}, it was observed that when the image of a coisotropic submanifold by a symplectic homeomorphism is a smooth submanifold, then this image is coisotropic and the symplectic homeomorphism intertwines the charateristic foliations of the source and target submanifolds. As a result, it acts on the reduction, which is locally a symplectic manifold, and one may ask which symplectic properties of the reduction are preserved. In other terms, the set of coisotropic submanifolds is $\cc^0$-rigid, as well as the characteristic foliations, and we ask which further symplectic invariants of coisotropic submanifolds are $\cc^0$-rigid. As far as we know, three results are already known, prior to this paper:
\begin{itemize}
\its When the reduction homeomorphism is smooth, it is symplectic \cite{buop}. 
\its Let $\Sigma,\Sigma'$ be hypersurfaces in symplectic manifolds, $h$ a symplectic homeomorphism that takes $\Sigma$ to $\Sigma'$ and $\red h$ its action at the level of the reduction. Then $\red h$ preserves the stable displacement energy of subsets of $\red \Sigma$ \cite{buop}. 
\its Let $h$ be a symplectic homeomorphism of $\T^{2n}$ that preserves the split coisotropic subtorus $\T^{2k+l}\times \{0\}$. Then $\red h:\T^{2k}\to \T^{2k}$ preserves the spectral capacity of open subsets of $\T^{2k}$ \cite{hulese2}.  
\end{itemize}
In these situations, the symplectic homeomorphism therefore also inherits symplectic properties at the level of the reduction. 
In this paper, we study the (non)-squeezing properties of the reduction of symplectic homeomorphisms in this setting, namely:
\begin{conj}\label{conj:rigred}
Let $h$ be a symplectic homeomorphism that takes some coisotropic submanifold to a smooth, hence coisotropic, submanifold. Then the reduction of $h$ does not squeeze balls into cylinders of smaller capacities. 
\end{conj}

Our first result generalizes the last result cited above. 
\begin{thm}\label{thm:c0rigredclosed} Let $N$ be a closed manifold that admits a metric with non-negative scalar curvature. Assume the existence of a symplectic homeomorphism defined in a \nbd of $B^{2k}(a)\times N\subset \C^k\times T^*N$ such that $h(B^{2k}(a)\times N)\subset Z^{2k}(A)\times N$. Then $A\geq a$. 
\end{thm}
In other terms, when the characteristic foliation consists of closed leaves and is trivial, and provided the leaves satisfy some topological property made explicit in remark \ref{rk:c0rigredclosed}, the reduction of a symplectic homeomorphism that preserves this submanifold verifies the non-squeezing property. 
Theorem \ref{thm:c0rigredclosed} relies on the classical non-squeezing theorem {\it via}  a short and easy argument. 

By contrast, our study of the general case is less conclusive. In a model situation, we still get the following result: 

\begin{thm}\label{thm:rigredloc}
Let $h$ be a symplectic homeomorphism defined in a \nbd of $B^{2k}(1)\times [-1,1]^{n-k}\subset \C^n$, with values in $\C^n$,
that takes the coisotropic submanifold $B^{2k}(1)\times [-1,1]^{n-k}$ to $\C^k\times \R^{n-k}$. Then there exists $\delta(h)>0$  such that for all $a<\delta(h)$, if $h\big(B^{2k}(a)\times [-1,1]^{n-k}\big)\subset Z^{2k}(A)\times \R^{n-k}$, $a\leq A$.   
\end{thm}
The dependence of the constant $\delta(h)$ with $h$ prevents the previous theorem to provide a satisfactory answer to conjecture \ref{conj:rigred}. The best we would get in the general framework of conjecture \ref{conj:rigred}, within the techniques of the present paper, would be that the balls of small diameters in the reduction of the source are either not-squeezed or highly distorted.  When the reduction has dimension $2$ however, non-squeezing means area-preservation,  and theorem \ref{thm:rigredloc} is enough for a complete answer. 
\begin{thm}\label{thm:rigred2} Conjecture \ref{conj:rigred} holds when the dimension of the reduction is $2$. In other terms, if $C$ is a coisotropic submanifold of dimension $n+1$ in a symplectic manifold $M^{2n}$, and if some symplectic homeomorphism takes $C$ to a smooth submanifold $C'$, then $C'$ is coisotropic, $h$ conjugates the characteristic foliations of $C, C'$ and $\red h:\red C\to \red C'$ is area preserving. 
\end{thm}

  
Notice that for a coisotropic submanifold, having a reduction of dimension $2$ means  being of minimal dimension above the lagrangian case. Thus, this corollary is in some sense orthogonal to the previous local result,  obtained in \cite{buop}, that concerned hypersurfaces, since these have maximal dimension below the open case. 
Theorem \ref{thm:rigredloc} relies on a slightly stronger version of the following coisotropic Camel theorem (see also \cite{bustillo} for a very close statement):

\begin{thm} \label{thm:camel}Let  $\phi:B^{2n}(a)\times \R^l\to \R^{2n}$ be a smooth map such that:
\begin{itemize}
\item[i)] $\forall t\in \R^l$, $\phi(\cdot,t)$ is a symplectic embedding,
\item[ii)] $\phi(\cdot,t)=\id+\sum_{1}^l t_i\frac\partial{\partial y_i}$ for $|t|\gg1$,
\item[iii)] $\im \phi\cap \{y_1=\dots=y_l=0\}\subset   [-1,1]^{n-k}\times  Z^{2(n-l)}(A) $.
\end{itemize}
Then $A\geq a$. 
\end{thm} 
For $l=1$, this is due to Eliashberg \cite{eliashberg}, and the generalization above is fairly straightforward. The strenghtened version necessary for our application, theorem \ref{thm:loccamel2}, requires more work, but is still not sufficient to handle conjecture \ref{conj:rigred}. In fact, this conjecture would follow if we could localize the symplectic camel theorem in the following sense:
\begin{conj} Theorem \ref{thm:camel} holds when iii) is replaced by:
\begin{itemize}
\item[iii)'] $\phi$ is knotted with $\partial Z_{A,R}^l$, where 
$$
Z_{A,R}^l:=\left\{(x_1,\dots,x_n,y_1,\dots,y_n)\; \left| 
\begin{array}{l}
y_1=\dots=y_l=0,\\ \pi(x_n^2+y_n^2)<A,\\ \|(x_1,\dots,x_{n-1},y_{l+1},\dots,y_{n-1})\|_\infty<R
\end{array}\right.
\right\}.
$$ 
\end{itemize}
In other terms, $\phi$ provides a symplectic isotopy of the ball $B^{2n}(a)$ through the window $Z_{A,R}^l$, without any further assumption that the isotopy avoids the coisotropic wall $\R^l\times \R^{2(n-l)}$ away from this window. 
\end{conj}

\paragraph{Organization of the paper:} The short section \ref{sec:rigredclosed} proves theorem \ref{thm:c0rigredclosed} 
 and is independent from the rest of the paper. In section \ref{sec:loccamel}, we state and prove two versions of the coisotropic versions of the Camel theorem that we have stated above. We provide an extensive argument, although some good references are already available \cite{mctr,niederkruger}. The reason is that these references either insist on the dimension $4$, or on the pseudoconvexity of the Camel space, which are not relevant in our context. Finally, we prove theorems \ref{thm:rigredloc} and  \ref{thm:rigred2} in section \ref{sec:coisorig}.

\paragraph{Notation:}
Throughout this paper, we adopt the following notation:
 \begin{itemize}
 \its $B(a)$ denotes the ball of capacity $a$ in $\C^n$.  
 \its Given two sets $A\subset B$, $\op(A,B)$ stands for an arbitrary but fixed \nbd of $A$ in $B$. 
 \its A continuous map $\Phi:\R^l\times B^{2n}(a)\to \C^n$ is a parametric symplectic embedding if $\Phi(t,\cdot):B^{2n}(a)\to \C^n$ is a symplectic embedding for all $t\in \R^l$.  We then denote $\Phi_t:=\Phi(t,\cdot)$
\its We split $\C^n$ as $\C^n=\C^{k+l}=\C\times \C^{k-1}\times \R^l\oplus i\R^l$ and we use coordinates  $(z_1,z',x,y)$ on these four factors, respectively. 
\its We endow $\C^{k-1}\times \C^l$ with the $L_\infty$-norm : $\|(z',x,y)\|:=\max \{\|z'\|_\infty, \|x\|_\infty,\|y\|_\infty \}$ (and $\|z'\|_\infty:=\max\{\|\re z'\|_\infty, \|\im z'\|_\infty\}$).
\its  $Z_{A,R}^{k,l}:=D(A)\times\{(z',0,y)\,|\, \|(z',0,y)\|<R\} \subset \C^k\times i\R^l$. When $k,l$ are clear from the context, we omit the subscripts for easier readability. 
\its $\Gamma_{A,R}^{k,l}:=S^1(A)\times\{(z',0,y)\,|\, \|(z',0,y)\|<R\} \subset \C^k\times i\R^l$ (this is the "horizontal" part of $\partial Z_{A,R}^{k,l}$),
\its $F_{A,R,M}^{k,l}:=D(A+\pi M^2)\times\{(z',x,y)\,|\, R-M<\|(z',0,y)\|<R+M \text{ and } \|x\|<M\} \subset \C^k\times i\R^l$. This is an $M$-\nbd of the "vertical" part of $\partial Z_{A,R}^{k,l}$. As a result, notice that any point in $F^{k,l}_{A,R,M}$ is at euclidean distance no less than $M'$ of $\partial F^{k,l}_{A,R,M+M'}$.  
\end{itemize}
The following figure illustrates our notation. It is rather accurate when $k=1$ and $l=0$. When $k+l\geq 2$, the figure is hopefully usefull, but it is only partially representative. The main difference is that $F_{A,R,M}$ becomes connected. 

 \begin{figure}[h!]
\begin{center}
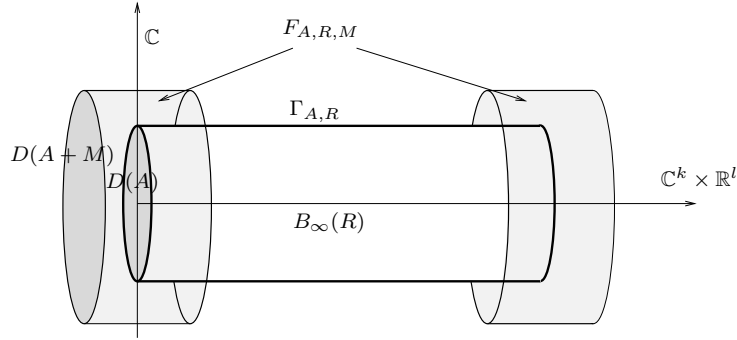
\end{center}
\caption{The sets $\Gamma_{A,R}$ and $F_{A,R,M}$ when $k=1$.}\label{fig:sft}
\end{figure}

\paragraph{Aknowledgements:}  Lev Buhovsky was initially engaged in this project and contributed to uncovering a crucial flaws in the initial argument. This paper therefore owes him a lot.

\section{Proof of theorem \ref{thm:c0rigredclosed}}\label{sec:rigredclosed}
Let $g$ be a metric on $N$ with non-positive sectional curvature,  $\eps>0$, and $h:\op(B^{2k}(a)\times N,\C^k\times T^*N)\to \C^k\times T^*N$ be a symplectic homeomorphism such that $h(B^{2k}(a)\times N)\subset Z^{2k}(A)\times N$. Consider an $\eps$-approximation of $h$ in the $\cc^0$-norm by a symplectic diffeomorphism $f$. Thus, $f$ is defined  in a \nbd of $B^{2k}(a)\times N$  in $B^{2k}(a)\times T^*N$ and provided this \nbd is small enough and the approximation fine enough, we have 
 $$
 f(B^{2k}(a)\times T_\delta^*N)\subset Z(A+\eps)\times T_{2\eps}^*N,
 $$ 
 where $\delta, \eps$ are small positive numbers. 
 Since $g$ has non-positive curvature, Cartan-Hadamard's theorem guarantees that the universal cover $\tilde N$ of $N$ is diffeomorphic to $\R^l$. The covering $\pi:\tilde N\to N$  lifts to symplectic coverings $T^*\pi:T^*\tilde N\to T^*N$ and $\Pi:=\id\times T^*\pi:\C^k\times T^*\tilde N\to \C^k\times T^*N$. The map $\pi$ obviously being a Riemanian covering with respect to the metrics $(\tilde g:=\pi^*g,g)$ on the pair $(\tilde N,N)$, it is easy to see that $\Pi^{-1}(U\times T_\delta^*N)=U\times T_\delta^*\tilde N$, where $T^*_\delta\tilde N$ stands for the $\delta$-\nbd of the zero section with respect to the lifted metric $\tilde g$.   Thus, $f$ lifts through $\Pi$ to a symplectic map 
$$
\tilde f:B(a)\times T^*_\delta \tilde N\lra Z(A+\eps)\times T^*\tilde N. 
$$ 
Since $h$ is a homeomorphism that sends $\{0\}\times N$ to $\{*\}\times N$, the projection $\text{pr}\circ h_{|{\{0\}\times N}}:\{0\}\times N\to N$ has degree $1$. This property also holds for $f$, and as a result, $\tilde f$ is injective, hence an embedding
$$
\tilde f:B^{2k}(a)\times T^*_{\delta} \tilde N\lra Z^{2k}(A+\eps)\times T^*\tilde N\approx Z^{2k}(A+\eps)\times T^*\R^l=Z^{2n}(A+\eps). 
$$

\begin{claim}\label{claim:negcurvsymp} 
$T^*_\delta\wdt N$ contains a symplectic ball of capacity $a$. 
\end{claim}

Let thus  $\phi:B^{2l}(a)\hra T^*_\delta \R^l$ be  a symplectic embedding. The map 
$\tilde f\circ (\id\times \phi)$ therefore provides a symplectic embedding of $B^{2k}(a)\times B^{2l}(a)$ into $Z(A+\eps)$. By Gromov's non-squeezing theorem, we  get $A+\eps\geq a$ and, letting $\eps$ to $0$, $A\geq a$. \cqfd

\noindent{\it Proof of claim \ref{claim:negcurvsymp}:}
Let $(N,g)$ be a non-positively curved Riemannian manifold and $(\wdt N,\tilde g)$ its universal cover. Fix any point $p\in \wdt N$.  By Cartan-Hadamard's theorem, the map  $f:=\exp_p:T_p\wdt N\to \wdt N$ is a diffeomorphism. Moreover, the negative curvature implies that the map
$$
D_\lambda(x):=f(\lambda f^{-1}(x)) 
$$
verifies $\|(D_\lambda')^{-1}\|^{-1}_{\tilde g}\geq \sqrt\lambda$. Indeed, given $u\in T_x\tilde N$ of norm $1$, $U(t):=t\mapsto D_t'(x)u$ is the Jacobi vector field along the geodesic issued from $p$ and passing at $x$ at time $1$, that vanishes at $0$ and equals $u$ at $1$. A classical computation shows that the function $h(t):=g(U(t),U(t))$ is convex in non-positive curvature.  Since it verifies $h(0)=0$ and $h(1)=1$, $h(t)\geq t$ for $t\geq 1$, so $\|D_\lambda'(x)u\|^2\geq \lambda$, which is our assertion. Given now any small ball $\cd\subset \wdt N$, $T^*\cd$ contains a relatively compact symplectic ball $\phi:B(a)\hra T^*\cd$. The image of $\phi$ lies in $T^*_{\tilde g,K}\wdt N$ for some large $K$. Putting $\lambda:=\frac \delta K$, the natural lift of $D_\lambda$ to $T^*\wdt N$ then takes $\phi(B(a))$ into $T^*_{\delta} \tilde N$. \cqfd

\begin{remark}\label{rk:c0rigredclosed}
The proof above shows that theorem \ref{thm:c0rigredclosed} holds whenever $N$ has the property that any \nbd $T^*_{\delta,\tilde g}\tilde N$ of the universal cover of $N$ equipped with some pull-back metric has infinite Gromov capacity. 
It is however slightly unclear to us whether there is a more conventional characterization of this natural symplectic property. The existence of a non-positive curvature on $N$ is enough, as explained above.  
\end{remark}

%

 \section{Two Camel theorems}\label{sec:loccamel}
 The aim of this section is to state and prove two Camel-like theorems, that we use in our proof of theorem \ref{thm:rigredloc}.
 \begin{theorem}\label{thm:loccamel2} There exists a universal constant $C>0$ for which the following holds. 
 Let $\Phi:\R^l\times B^{2n}(a)\to \C^n=\C^{k+l}$ be a parametric symplectic embedding which verifies:
 \begin{itemize}
  \item[(i)] $\Phi$ is standard at infinity: $\exists K>0$ such that for $|t|>K$, $\Phi_t(z)=z+t\frac\partial{\partial t}$, where
 we see $t\in \R^l\simeq \{0\}\times \R^l\subset \C^k\times\C^l$,\vspace*{-,2cm}
 \item[(ii)] $\Phi$ is knotted with $\partial Z_{A,R}^{k,l}$. 
 \item[(iii)] $\im \Phi$ does not intersect $F_{A,R,CA}^{k,l}$. 
 \end{itemize} 
 Then $A\geqslant a$. 
 \end{theorem}
 
 \begin{theorem}
 \label{thm:loclagcamel} There exists a universal constant $C>0$ for which the following holds. 
 Let $L$ be a Lagrangian embedding of $(S^1)^k\times \R^{l}$ into $\C^n=\C^{k+l}$ which verifies:\vspace*{-,2cm}
\begin{itemize}
\item[(i)] $L$ is properly embedded, standard at infinity: there is a compact set $K\subset \C^n$ such that \vspace*{-,3cm}
$$
L\cap {}^c K=\big(S^1(a)^k\times \R^{l}\big) \cap {}^c K =\{\pi|z_i|^2=a, i=1\dots k\}\times \R^l\cap {}^c K \subset \C^k\times \C^l.\vspace*{-,2cm}
$$
\item[(ii)] $L$ is knotted with $\partial Z_{A,R}^{k,l}$.
\item[(iii)] $L$ does not meet $F_{A,R,CA}^{k,l}$. 
\end{itemize}
Then $a\leqslant A$.
\end{theorem}

Obivously, the assumptions on theorem \ref{thm:loccamel2} are less restrictive than the usual ones, so theorem \ref{thm:loccamel2} implies theorem \ref{thm:camel}. Theorem \ref{thm:loclagcamel} is also a slightly localized version of an existing Lagrangian camel theorem \cite{bustillo}. As already explained, getting a genuinely localized version of both results would consist in removing assumption (iii).  
 We prove both theorems below, because both seem interesting in their own sake, and that none of them seems to follow from the other, but their proofs are very similar. For our application to the $\cc^0$-rigidity of the reduction of a symplectic homeomorphism however, only theorem \ref{thm:loccamel2} is relevant, through the following corollary:
 \begin{corollary}\label{cor:camelcoiso}
If there is a compactly supported symplectic diffeomorphism $f$ of $\C^n=\C^{k+l}$ such that 
\begin{itemize}
\item[(i)] $f(B^{2k}(a)\times \R^l)$ is knotted with $\partial Z_{A,R}^{k,l}$ for some $R$,
\item[(ii)] $f(B^{2k}(a)\times \R^l)\cap F_{A,R,CA}^{k,l}=\emptyset$, 
\end{itemize}
 then $A\geq a$.     
 \end{corollary}
\noindent{\it Proof:} Let $K,\eps>0$ to be  taken very large and very small, respectively. Let $\rho:\R^+\to [\eps,1]$ be a smooth function with $\rho(s)\equiv \eps$ for $s<K$ and $\rho(s)\equiv 1$ for $s>K+1$. The map \fonction{\phi}{B^{2k}(a)\times B^{2l}(a)\times \R^l}{\C^{k+l}}{(z_1,z_2,t)}{(z_1,\frac 1{\rho(\|t\|)}\re z_2+t,\rho(\| t\|)\im z_2)}
verifies:
\begin{itemize}
\item[(i)] $\phi$ is smooth,\vspace*{-,2cm}
\item[(ii)] for all $t\in \R^l$, $\phi(\cdot,t)$ is a symplectic embedding of $B^{2k}(a)\times B^{2l}(a)\supset B^{2n}(a)$,\vspace*{-,2cm}
\item[(iii)] for $|t|>K+1$, $\phi(z_1,z_2,t)=(z_1,z_2)+t$,\vspace*{-,2cm}
\item[(iv)] for $|t|<K$, $\phi(z_1,z_2,t)$ belongs to an $\eps$-\nbd of $B^{2k}(a)\times \R^l\subset \C^k\times \C^l$. \vspace*{-,2cm}
\end{itemize}
Let now $f$ be the symplectic diffeomorphism considered by our statement. For $\eps\ll 1$ and $K\gg 1$, the restriction of the map $\Phi:=f\circ \phi$ to $B^{2n}(a)\times \R^l\subset B^{2k(a)}\times B^{2l}(a)\times \R^l$ verifies all the assumptions of theorem \ref{thm:loccamel2}, so  $A\geq a$.\cqfd


Before proving theorems \ref{thm:loccamel2} and \ref{thm:loclagcamel}, let us unify their statements. First, we fix $k,l,n$ so we drop all these indices from our notations. Notice also that scaling the symplectic form allows to assume that $A=1$, which we do henceforth, removing from our notations any reference to $A$ (hence $Z_R,\; F_{R,M},\Gamma_{R}$ stand for $Z_{1,R}^{k,l},F_{1,R,M}^{k,l},\Gamma_{1,R}^{k,l}$). Let $f:X\times \R^l\to \C^n=\C^k\times \C^l$ be a map that is either a Lagrangian embedding if $X=S^1(a)^k$, or a parametric symplectic embedding of the ball if $X=B^{2n}(a)$. It is also standard at infinity, knotted with  $\partial Z_R$, and its image avoids $F_{R,M}$.  We need to find a constant $C$ such that $M\geq C$ implies the inequality $a\leq 1$. 
Since $f$ is standard at infinity in both cases, its image is contained in  $\{(z,x,y)\;|\;\| (z,y)\|<K\}$ for some large $K$, that we take larger than $M$. 
Moreover, still because $f$ is standard at infinity,  we can compactify the triple $(f,\dom f,\C^k\times i\R^l\oplus \R^l)$ by compactifying the $\R^l$-factor to a very large torus, which we choose in any case to contain $\{|x|\geq \max\{1,M\}\}$. We get a new triple $(\tilde f,\dom \tilde f,\C^k\times i\R^l\times \T^l)$, where $\dom \tilde f=X\times \T^l$, $X=S^1(a)^k$ or $X=B^{2n}(a)$, and $\tilde f$ is either a Lagrangian embedding or a parametric symplectic embedding. Moreover, $\im \tilde f$ is still knotted with $\partial Z_R$, and avoids $F_{R,M}$, both sets being unaffected by the compactification because $Z_R\subset \{x=0\}$ and $K>M$.  In order to keep light notation, we drop the tilde from $\tilde f$. We call $Y:=\C^k\times i\R^l\times \T^l$ our ambient space. It is equipped with the standard symplectic form $\om$, and with a standard complex structure $J_\st$. The almost complex structures that we will consider in this analysis belong to the set of $\om$-compatible almost complex structure $\cj(\om)$ on $Y$ and are of the following type. First, fix two disjoint \nbds $V_\Lambda\Subset \{\|(z,y)\|<K\}$ of $\Lambda:=\im (f)$ and $V_\Gamma$ of $\Gamma_{R}\cup F_{R,M}$, and consider some $\om$-compatible almost complex structure $J_\Lambda$ on $V_\Lambda$. 
We define 
$$
\cj(R,M,J_\Lambda):=\{J\in \cj(\om)\,|\, J_{|V_\Lambda}\equiv J_\Lambda\text{ and }J_{|V_\Gamma\cup \{\|(z,y)\|>K\}}\equiv J_\st\}.
$$ 

The idea is the same as for the classical Camel theorem: we must fill $\Gamma_{R}$ by holomorphic discs of area $1$. Since $\Gamma_{R}$ is not totally real, the problem is not in a good setting. We therefore introduce for $c\in \R^{k-1}$
$$
\Gamma_{R,c}:=\Gamma_{R}\cap \{\re z'=c\}=S^1\times \{(z',0,y) \; |\; \re z'=c,\; \|(\im z',y)\|<R\}
$$
which provides a foliation of $\Gamma_R$ by Lagrangian leaves, which we will each fill. 
The main difficulty here is that  we do not assume that the map $f$ avoids the coisotropic wall $\C^k\times i\R^l$ away from $Z(1)\times i\R^l$, so the classical filling technique must be adapted. For the reader aware of the classical proof, it might be worth having in mind that the main issue in this setting, compared to the more classical one, is the compactness. Since $J=J_\st$ on $F_{R,M}$, the circles $S^1\times \{(z',0,y)\}$ are filled by the vertical discs $\D\times\{(z',0,y)\}$ for $\|(z',0,y)\|\approx R$. We need to guarantee that the boundaries of the non-vertical holomorphic discs do not approach $\partial \Gamma_{R,c}$, which will be realized by taking $M$ large enough  (see lemma \ref{le:trans1}). 

Although the proof of the Camel theorem is folklore, we could not find a complete proof in the litterature, except in dimension $4$ \cite{mctr}. We therefore provide a complete proof, and not only of the compactness issue.

\subsection{Preliminaries for the compactness}\label{sec:blowup}
Here is the main statement we will need in order to address the compactness issues. 
\begin{prop}\label{prop:maincompact} 
There exists $M>0$ such that for any $R$, $c$, $J_\Lambda\in \cj(\om)$ and $J\in \cj(R,4M,J_\Lambda)$, all the $J$-holomorphic discs with boundary on $\Gamma_{R,c}$, of area $1$, that intersect $F_{R,M}$ are vertical discs. 
\end{prop}
This statement  relies on classical area estimation of analytic sets, provided by the next  lemmas. 
\begin{lemma}\label{le:goodball} There exists a relative holomorphic embedding $\phi:(B^{2n}(\eps_0),\R^n)\hra  (Y,\Gamma_0)$, with $\phi(0)=(z_1=1,z'=x=y=0)\in \Gamma_{0}$. 
\end{lemma}
\noindent {\it Proof:} There exists a holomorphic embedding of $(\D,\R)\to (\C,S^1)$ because $S^1$ is real-analytic. Taking product with the identity  provides the desired holomorphic map. \cqfd

In the next lemma, we say that an analytic subset has {\it real boundary} if every local branch of $X$ near a point $x\in \partial X$ can be parametrized by a holomorphic map defined on $\D\cap\{\im z\geq 0\}$. For us, the main example of an analytic subset with real boundary on $\Gamma_c$ is the image of a holomorphic disc $u:(\D,\partial \D)\to (Y,\Gamma_c)$. 

\begin{lemma}\label{le:monotonicity} 
Let $\om$ be a symplectic form on $B^{2n}(1)\subset \C^n$, compatible with $J_\st$. There exists a constant $\hbar$, that depends only on $\om$,  such that:
\begin{enumerate}
\item[(i)] For any proper analytic subset $X$ of $B^{2n}(1)$ of complex dimension $1$,  with $0\in X$, $\ca_\om(X)\geq \hbar$.
\item[(ii)] For any proper analytic subset $X$ of $B^{2n}(1)$ of complex dimension $1$, with real boundary on $\R^n$ and $0\in \partial X$, 
$\ca_\om(X)\geq \frac \hbar 2$. 
\end{enumerate}
\end{lemma}
 \noindent {\it Proof:} Both assertions are well-known when $\om=\om_\st$ (this is the monotonicity lemma). Since $\om_\st$ and $\om_0:=\om(0)$ are $J_\st$-compatible, they have symplectic basis of the form $(e_1,J_\st e_1,\dots,e_n,J_\st e_n)$, $(f_1,J_\st f_1,\dots,f_n,J_\st f_n)$. The $J_\st$-linear map $A$ that takes $e_i$ to $f_i$ therefore transports $\om_0$ to $\om_\st$. Then $A(X)$ is an analytic subset passing through $0$, proper in $A(B(1))\supset B(r)$, for some $r>0$ that depends only on $A$, hence on $\om_0$. 
 
 Let us prove (i). The montonicity lemma guarantees that $\ca_{\om_\st}(AX\cap B(\eps))\geq \pi \eps^2$ for all $\eps<r$. Since moreover $A^*\om_0=\om_\st$, we have $A^*\om=\om_\st +R$, where $R\in O(\eps)$ on $B(\eps)$, so 
 $$
 \ca_{A^*\om}(AX\cap B(\eps))=\int_{AX\cap B(\eps)}A^*\om=\int_{AX\cap B(\eps)} \om_\st+R\geq\ca_{\om_\st}(AX\cap B(\eps))-C\eps\ca_{g_\st}(AX\cap B(\eps)),
 $$
where $\ca_{g_\st}$ stands for the euclidean area, and $C$ depends only on $\om$. Now since $AX$ is an analytic set, $\ca_{g_\st}(AX\cap B_\eps)=\ca_{\om_\st}(AX\cap B_\eps)$, so we get 
$$
\ca_{A^*\om}(AX\cap B(\eps))\geq \ca_{\om_\st}\big(AX\cap B(\eps)\big)(1-C\eps)\geq \pi\eps^2(1-C\eps).
$$ For some small enough $\eps_0<r$, we therefore see that $\ca_{A^*\om}(AX\cap B(\eps_0))\geq \frac \pi 2\eps_0^2$. Finally, since $\eps_0<r$, $B(r)\subset A(B(1))$ and $X$ is analytic, 
$$
\ca_\om(X\cap B(1))\geq \ca_\om(X\cap A^{-1}(B(\eps_0)))=\ca_{A^*\om}(AX\cap B(\eps_0))\geq \frac \pi 2\eps_0^2=\hbar. 
$$
  
 In order to prove (ii), we consider a proper analytic set with real boundary on $\R^n$. Recall that this means that $(X,\partial X)\subset (B^{2n}(1),\R^n)$, and that  every local branch $B$ of $X$ near a point $x\in \partial X$ can be parametrized by a holomorphic function $f_B$ defined on a \nbd of $0$ in $\H$, and we assume here that $f_B(\R)\subset \R^n$. Then $X\cup \sigma(X)$ (where $\sigma(z_1,\dots,z_n)=(\bar z_1,\dots,\bar z_n)$) is analytic on $X$, on $\sigma(X)$, and around the points of $\partial X=\partial \sigma(X)$ by the reflection principle, applied to each component of $f_B$.  Thus $X\cup \sigma(X)$ is a proper analytic set of $B^{2n}(1)$, without boundary, which passes through $0$, so its area is at least $\hbar$ by the first part of the argument. The first part of the argument also shows that 
 $$
  \ca_\om(X)\geq (1-C\eps_0)\ca_{\om_\st}(X\cap B(\eps_0))=\frac 12(1-C\eps_0)\ca_{\om_\st}(X\cup \sigma(X)\cap B(\eps_0))\geq \frac 12\hbar.\hspace{1cm}\square  
  $$

\begin{lemma}\label{le:idiot} A non-constant $J_\st$-holomorphic disc with boundary on $\Gamma_c$ is a branched covering of a vertical disc, where  the vertical discs are understood to be the discs $\D\times \{(\re z'=c,\im z'=c_1,y=c_2)\}$.
 \end{lemma}
  \noindent{\it Proof:} We can write $u=(u_1,u',u_2)$, where $u_1:(\D,\partial \D)\to (\C,S^1)$,  $u':(\D,\partial \D)\to (\C^{k-1},\{\re z'=c\})$, $u_2:(\D,\partial \D)\to (i\R^l\times \T^l, \{x=0\})$ are $J_\st$-holomorphic discs. Then $u',u_2$ are easily seen to be constant discs by Stokes formula for instance (the integrals of $\lambda_\st$ over their boundaries vanish). And $u_1$ takes values in $\D$ by the maximum principle. Finally,  $u_1:\D\to \D$ is a proper map, so is a branched covering of $\D$.\cqfd
  
%
  
 \noindent{\it Proof of proposition \ref{prop:maincompact}:} Let $M,R>0$, $J_\Lambda\in\cj(\om)$, $J\in \cj(R,4M,J_\Lambda)$ and $u:(\D,\partial \D)\to (Y,\Gamma_{R,c})$ be a $J$-holomorphic disc of area $1$ that intersects $F_{R,M}$. Notice that if this disc lies in $F_{R,4M}$, it is $J_\st$-holomorphic hence vertical by lemma \ref{le:idiot}. We therefore assume throughout this proof that $u(\D)$ meets both $F_{R,M}$ and $Y\priv F_{R,4M}$. We are then in one of three possible situations:
 \begin{itemize}
 \its Either $u(\partial \D)\subset F_{R,3M}$. Then since $\im u$ meets $Y\priv F_{R,4M}$, there is a euclidean ball $B$ of radius $\nf M2$ centered on $\im u$ that lies in $F_{R,4M}\priv F_{R,3M}$. Then $u(\D)\cap B$ is a proper analytic subset of $B$ (because $J=J_\st$ on $F_{R,4M}$), so has area at least $\nf{\pi M^2}4$ by the monotonicity lemma. Since $u$ has total area $1$, this situation does not occur if we choose $M\geq 2$. 
 \its Or $u(\partial \D)\cap F_{R,2M}=\emptyset$. Since $u$ visits $F_{R,M}$, the same argument shows that $u$ has area at least $\nf {\pi M^2}4$, so does not occur again if $M\geq 2$. 
 \its Or $u(\partial \D)$ meets both $F_{R,2M}$ and $Y\priv F_{R,3M}$. Denote by $\delta_0$ the euclidean diameter of the $J_\st$-holomorphic ball $\phi:(B(\eps_0),\R^n)\to (Y,\Gamma_0)$ centered at $p_0:=(z_1=1,x=y=z'=0)$ provided by lemma \ref{le:goodball}. By assumption, $u(\partial \D)$ is connected and meets both $F_{R,2M}$ and $Y\priv F_{R,3M}$ so its intersection with $F_{R,3M}\priv F_{R,2M}$ has diameter at least $M$. Thus, we can center at least $k:=\lfloor \nf M {2\delta_0}\rfloor$ disjoint euclidean balls of radii $\delta_0$ on points $p_j\in u(\partial \D)\cap F_{R,3M}\priv F_{R,2M}$ ($j\in[1,k]$), and these balls lie in $F_{R,4M}$ if $M\geq \delta_0$.  Denote by $\tau_j$ the composition of a translation  in the $(z',y)$-factor and rotation in the $z_1$-factor that brings $p_0$ to $p_j$. Then, the $\tau_j\circ \phi(B(\eps_0))$ provide $k$ disjoint $J_\st$-balls centered on the $p_j$. Moreover, since $\om$ is invariant by the $\tau_j$,  the pull-backs $(\tau_j\circ \phi)^*\om$ are a fixed symplectic form on $B^{2n}(\eps_0)$. As a result, the area of $u$ is at least the sum of the areas of the intersections of $u(\D)$ with these $k$ balls, each of which being at least some constant $\hbar$ by lemma \ref{le:monotonicity}.(ii) (this constant depends only on $\phi$ by the discussion above). Thus the area of $u$ is at least $k\hbar\geq \big(\nf{M}{2\delta_0}-1\big)\hbar$. Taking $M> 2\delta_0\big(1+\nf 1\hbar\big)$, the area of $u$ exceeds $1$, so this situation does not happen neither.   
  \end{itemize} 
 We finally conclude that if $M> \max\{2,\delta_0,2\delta_0\big(1+\frac 1\hbar\big)\}$, the only discs of area $1$ with boundary on $\Gamma_{R,c}$ that intersect $F_{R,M}$ are the vertical ones. \cqfd

\subsection{Filling by holomorphic discs}  
Recall our notations: $\cj(\om)$ is the set of $\om$-compatible almost complex structures on $Y=\C^k\times i\R^l\times \T^l$, $Z_R=\D\times \{(z',0,y),\; \|(z',y\|\leq R)\}$, $\Gamma_R=S^1\times\{(z',0,y),\; \|(z',y\|\leq R)\} $ is the "horizontal" part of $\partial Z_R$, while $F_{R,M}=\D\times \{(z',x,y),\; \|x\|\leq M,\; R-M\leq\|(z',y)\|\leq R+M\}$ is the $M$-\nbd of the vertical part of $\partial Z_R$. The constant $M$ is chosen large enough so that the conclusion of proposition \ref{prop:maincompact} holds, and $R$ is then chosen large compared to $M$, so that $F_{R,M}$ does not cover $Z_R$.   A map $f:X\times \T^l\to Y$ is given, whose image is surrounded by a small \nbd $V_\Lambda$. We consider $M$ for which proposition \ref{prop:maincompact} holds, and we assume that $V_\Lambda\subset Y\priv (\Gamma_R\cup F_{R,4M})$.   We aim at filling $\Gamma_R$ with holomorphic discs for good complex structures in $\cj(\om)$.  In order to apply the Fredholm theory, we rather fill each Lagrangian leaf $\Gamma_{R,c}$, $c\in \cd:={\bar B}^{k-1}(R)\subset \R^{k-1}$ in a consistent way. To prove theorem \ref{thm:loccamel2}, we will further need a parametric version of this filling. Throughout this section, $T$ stands for  a parameter space, which is a smooth closed manifold, $J_\Lambda:T\times \cd\to \cj(\om)$ is any $\cc^l$-smooth map, and $J:T\times \cd\to \cj(\om)$ is a smooth map that verifies $J_{(c,t)|V_\Lambda}\equiv J_\Lambda$ and $J\equiv J_\st$ on $F_{R,4M}$.  Formally, $J$ is a $\cc^\ell$-section of the bundle $\cj(R,4M,J_\Lambda)\to T\times \cd$, whose fiber is 
$$
\cj(R,4M,J_\Lambda(t,c)):=\{J\in \cj(\om),\; J_{|V_\Gamma\cup\{\|(z,y)\|>K\}}\equiv J_\st,\; J_{|V_\Lambda}\equiv J_\Lambda(t,c)\},
$$
where $V_\Gamma$ is any fixed \nbd of $\Gamma\cup F_{R,4M}$. 
 Also, $\pi_2(Y,\Gamma_c)$ is generated by the class of vertical discs (parametrized injectively and holomorphically). The image of this class by the Hurewitz morphism $\pi_2(Y,\Gamma_c)\to H_2(Y,\Gamma_c)$ is denoted by $E_c$. The spaces of interest for us are
$$
\cm(J):=\{(t,c,u),\; (t,c)\in T\times \cd,\, u:(\D,\partial \D)\to (Y,\Gamma_{R,c}), \, \bar \partial_{J_{(t,c)}}u=0,\, [u]=E_c\}
$$
and what will turn out to be the interior of $\cm(J)$, obtained by replacing the closed ball $\cd$ by the open ball $\rond \cd:=B^{k-1}(R)$ and $\Gamma_{R,c}=S^1\times\{(\im z',c,y),\; \|(\im z',y)\|\leq R\}$ by $\rond \Gamma_{R,c}:=S^1\times\{(\im z',c,y),\; \|(\im z',y)\|< R\}$:
$$
\cm'(J):=\{(t,c,u),\; (t,c)\in T\times \rond \cd,\, u:(\D,\partial \D)\to (Y,\rond\Gamma_{R,c}), \, \bar \partial_{J_{(t,c)}}u=0,\, [u]=E_c\}
$$
Denoting $\tau_c$ the translation of vector $c\in \R^{k-1}$, the spaces $\cm(J)$ and $\cm'(J)$ are in one-to-one correspondence via the map $(t,c,u)\to (t,c,\tau_{-c}\circ u)$ with the following spaces, more suited to our analysis:
$$
\begin{array}{l}
\cm_0(J):=\{(t,c,u),\; (t,c)\in T\times \cd,\, u:(\D,\partial \D)\to (Y, \Gamma_{R,0}),  \, \bar \partial_{\tau_c^*J_{(t,c)}}u=0,\, [u]=E_0\}.\\
\cm_0'(J):=\{(t,c,u),\; (t,c)\in T\times \rond \cd,\, u:(\D,\partial \D)\to (Y,\rond \Gamma_{R,0}),  \, \bar \partial_{\tau_c^*J_{(t,c)}}u=0,\, [u]=E_0\}.
\end{array}
$$  
 Notice  that the almost complex structures $\tau_c^*J_{(t,c)}$ do not belong anymore to $\cj(R,4M,J_\Lambda(t,c))$ but to 
 $$
 \cj(R,4M,c,J_\Lambda):=\{J\in \cj(\om)\;| \; J_{|V_\Lambda}\equiv \tau_c^*J_\Lambda(t,c)), J_{|\tau_c^*(V_\Gamma\cup\{\|(y,z)\|>K\})}\equiv J_\st\}.
 $$
 Also, $\cm(J)$ is endowed with an action of $G:=\psl_2(\R)\simeq \aut(\D,j)$ by source reparametrization of $u$. 
  
 \begin{lemma}\label{le:trans1}
 If $(t,c,u)\in \cm_0(J)$ and $\im u\subset \tau_c^{-1}(F_{4R,M})$, $u$ is a vertical disc, and $\tau_c^*J_{(t,c)}$ is Fredholm-regular with respect to $u$. 
 \end{lemma}
 The Fredhlom theory and its usual notations are assumed in the following proof. The reader can consult \cite{mcsa2}, or look at the proof of lemma \ref{le:smoothmoduli} for few details. 
 
 \noindent{\it Proof:} $\tau_c^*J_{(t,c)}=J_\st$ on $\tau_c^{-1}(F_{4R,M})$, so the first part of the assertion directly follows from lemma \ref{le:idiot}. 
Moreover $\tau_c^*J_{(t,c)}$ is linear (constant) on this region, so the derivative 
$$
D_u\xi:=D\bar\partial_{\tau_c^*J_{(t,c)}}(u)\xi=\bar\partial_{\tau_c^*J_{(t,c)}}\xi=\bar\partial\xi.
$$
Therefore, if $\xi$ belongs to $\ker D_u$,  $\xi$ is easily seen to be a constant in an $n-1$-real parameter space, once modded out the source reparametrization. It follows that $\dim \ker D\bar \partial_{\tau_c^*J_{(t,c)}}(u)=n-1$, while the index of $D\bar \partial_{\tau_c^*J_{(t,c)}}(u)$ is $n-1$, so $D\bar \partial_{\tau_c^*J_{(t,c)}}(u)$ is indeed surjective.\cqfd

\begin{lemma}\label{le:sominj}
If $(t,c,u)\in \cm_0(J)$, then $u$ is somewhere injective, and almost all points of $\im u$ have exactly one preimage. 
\end{lemma}
This point is pefectly clear because the class $E$ has least area in $H_2(M,\Gamma_c)$. Here is nonetheless a full proof.

\noindent{\it Proof:}   We first recall that the critical set ${\cal C}(u)$ of a non-constant $J$-holomorphic disc $u$ is a discrete subset of $\D$, hence negligible (see \cite[Lemma 2.4.1]{mcsa2}).  
Let $(t,c,\tilde u)\in \cm_0(J)$, so $ u:=\tau_c\circ \tilde u$ is $J_{(t,c)}$-holomorphic. By \cite{lazzarini}, there exists a holomorphic disc $v:(\D,\partial \D)\to (Y, \Gamma_c)$ with $\im v\subset \im u$,  and $v^{-1}(v(z))=\{z\}$ for almost all points $z\in \D$. Then, 
$$
0<\ca_\om(v)\leq \ca_\om(u)=[\om][E]=1.
$$
Since moreover $\hur_2(Y,\Gamma_c)$ is generated by the class $[E]$, it follows that $\ca_\om(v)=A=\ca_\om(u)$. Let $D':=u^{-1}(v(\D))$. This is  an open subset of $\D$ because each local branch of $u^{-1}\circ v$ is holomorphic (its derivatives are $\C$-linear maps) and non-constant. We also  claim that it has full measure in $\D$. Indeed, 
$$
1=\ca_\om(u)=\int_{D'}u^*\om+\int_{\D\priv D'}u^*\om\geqslant \ca_\om(v)+\int_{\D\priv D'}u^*\om\geq \ca_\om(v)=1.
$$
This chain of inequalities consists therefore of equalities only, so 
$$
0=\int_{\D\priv D'}u^*\om=\int_{\D\priv D'\cup \, {\cal C}(u)}u^*\om.
$$
On $\D\priv(D'\cup {\cal C}(u))$, $u$ is a local diffeomorphism and $u^*\om>0$. Thus 
$$
0=\leb\big((\D\priv D')\priv {\cal C}(u)\big)\geq \leb(\D\priv D')-\leb ({\cal C}(u))=\leb(\D\priv D'),
$$ 
which shows that $D'$ has full measure. Define now 
$$
\textnormal{NI}:=\{p\in v(\D)\, |\, \#v^{-1}(\{p\})\neq 1\}, \hspace{,3cm}N_v:=v^{-1}(\overline{\textnormal{NI}}),\hspace{,3cm} N_u:=u^{-1}(\overline{\textnormal{NI}}).
$$
 It is clear that $\overline{\textnormal{NI}}\subset \textnormal{NI}\cup v({\cal C}(v))=v(v^{-1}(\textnormal{NI})\cup {\cal C}(v))$. Since ${\cal C}(v)$ and $v^{-1}(\textnormal{NI})$ are both negligible sets, and $v$ is smooth, it follows that $\overline{\textnormal{NI}}$ 
 has vanishing $2$-dimensional Hausdorff measure. Then 
$$
N_u=u^{-1}(\overline{\textnormal{NI}})=\left(u^{-1}(\overline{\textnormal{NI}})\cap {\cal C}(u)\right)\cup \left(u^{-1}\big(\overline{\textnormal{NI}}\priv u({\cal C}(u))\big)\right). 
$$
Since ${\cal C}(u)$ is discrete, the measure of the first set in the right hand side vanishes. On the other hand, $u^{-1}$ is smooth on $\overline{\textnormal{NI}}\priv u({\cal C}(u))$, so it preserves the vanishing of the $2$-dimensional Hausdorff measure, so the latter is also Lebesgue-negligible. We therefore conclude that $\leb(N_u)=0$.  
Now $\phi:=v^{-1}\circ u:D'\priv N_u\lra \D\priv N_v$ is a holomorphic map, surjective on $\D\priv N_v$ because $\im v\subset \im u$, so the map $\deg \phi:\D\priv N_v\to \Z$ that associates to each point of $\D\priv N_v$ the algebraic count of its preimages takes values in $\N^*$. Thus, 
$$
1=\int_\D u^*\om=\int_{D'\priv N_u}u^*\om=\int_{D'\priv N_u}(v\circ \phi)^*\om=\int_{\D\priv N_v}\deg \phi \;v^*\om\geq  \int_\D v^*\om=1
$$ 
(the last inequality holds because $v^*\om>0$). Thus $\deg \phi=1$ almost everywhere in $\D\priv N_v$, therefore constantly equals $1$ because $\phi$ is holomorphic. Thus $\phi=v^{-1}\circ u$ is injective on $D'\priv N_u$, so also is $u_{|D'\priv N_u}$. This proves the lemma because $D'\priv N_u$ has full measure and $u$ is smooth.\cqfd


 \begin{lemma}\label{le:smoothmoduli}
 For every smooth section $J\in \Gamma^\ell( \cj(R,M,J_\Lambda))$, and for every $\eps>0$, there exists  $J'\in \Gamma^\ell(\cj(R,M,J_\Lambda))$, with $d_{\cc^\ell}(J,J')<\eps$, such that $\cm_0'(J')$ is a smooth submanifold. 
 \end{lemma}
 \noindent{\it Proof :} This is standard and lengthy, so we only sketch the proof briefly, referring to \cite[\S 3]{mcsa2} for details. One defines the universal moduli space 
 $$
\cm:=\{(u,J,t,c),\; (t,c)\in T\times \rond \cd,\; J\in \cj(R,M,J_\Lambda),\, u:(\D,\partial \D)\to (Y,\rond\Gamma_{R,0}),\, \bar\partial_{\tau_c^*J_{(t,c)}}u=0,\,[u]=E_0\}. 
$$
The first point is to see that  $\cm$ is a Banach manifold. Let $\cb^{k,p}:=W^{k,p}((\D,\partial \D)\to (Y,\rond \Gamma_{R,0}))$, $\ce^{k-1,p}$ the Banach vector bundle over $\cb^{k,p}\times \cj(R,M,J_\Lambda)\times T\times \rond \cd$ whose fiber is 
$$
\ce^{k-1,p}(u,J,t,c):=W^{k-1,p}\big((\D,\partial \D)\to (\Lambda^{(0,1)}_{\tau_c^*J}(u^*TY),\Lambda^{(0,1)}_{\tau_c^*J}(u^*T\Gamma_{R,0}))\big).
$$
Then $\cf(u,J,t,c):=\bar \partial_{\tau_c^*J}u$ defines a differentiable section of this Banach bundle, whose zero-section is precisely $\cm$. The derivative of $\cf$ at the zero-section can be computed explicitely (using any connection on $\cj(R,M,J_\Lambda)$, we identify $T_{(t,c,J)}\cj(M,J_\Lambda)$ with a subset of $ T_tT\times T_c\R^{k-1}\times T_J\cj(\om)$):
$$
d\cf(\xi,Y,\delta t,\delta c)=D_u\xi+\frac 12\tau_c^*\left(Y(u)+\frac{\partial J}{\partial c}(u)\delta c +\frac{\partial J}{\partial t}(u)\delta t\right)\circ du\circ j,
$$
where $D_u(\xi):=D\bar \partial_{\tau_c^*J}(u)\xi$. 
We claim that this differential is surjective at all $(u,t,c,J)\in \cm$. 
  Since $D_u$ is onto for the elements $u\in \cm$ which are contained in $\tau_c^{-1}(F_{R,4M})$ by lemma \ref{le:trans1}, $d\cf(u,J,t,c)$ is onto for all these elements. Let $(u,J,t,c)$ be another element of $\cm$. To show that $d\cf(u,J,t,c)$ is surjective, consider an element $\eta$ of the dual that anihiliates the image of  $d\cf(u,J,t,c)$. Then in particular, for all elements $(\xi,Y)\in T_{(u,J)}\cb^{k,p}\times \cj(M,J_\Lambda(t,c))$, 
$$
\langle D_u\xi,\eta \rangle_{L^2}=0 \text{ and }\langle Y(u)\circ du\circ j,\eta \rangle_{L^2}=0. 
$$  
The first equation guarantees that if $\eta$ vanishes on some region, it vanishes identically. Now $\im u$ intersects the region $Y\priv(F_{R,4M}\cup V_\Lambda)$, where, by lemma \ref{le:sominj} there exists a somewhere injective point of $u$, and where $J$ is unconstrained by definition of $\cj(M,J_\Lambda(t,c))$. Using the freedom on $Y$ at such a point, the classical argument  shows that $\eta\equiv 0$. It follows that $d\cf(u,J,t,c)$ is surjective also at this point $(u,J,t,c)\in \cm$. Thus, $0$ is a regular value of $\cf$ and $\cm=\cf^{-1}(0)$ is a Banach manifold.

Now the map $\pi:\cm\to \cj(R,M,J_\Lambda)$ is a Fredholm map, so for any smooth section $J:T\times \rond \cd\to \cj(R,M,J_\Lambda)$, there is a $\cc^\ell$ $\eps$-small deformation $P'$ of the submanifold $P:=\{(J(t,c),t,c)\}$, such that $\pi$ is transverse to $P'$ \cite{smale2}. Such a perturbation is obviously the graph over $T\times \cd$ of a section $J'$, which verifies $d_{\cc^\ell}(J,J')<\eps$. The transversality ensures that $\pi^{-1}(P')=\cm_0'(J')$ is a smooth manifold. \cqfd
We will call henceforth {\it generic} the sections $J:T\times \cd\to \cj(R,M,J_\Lambda)$ such that $\cm_0'(J)$, hence $\cm'(J)$ is a smooth submanifold. We now address the compactness. 
\begin{lemma}\label{le:compactmoduli}
Let $J\in \cj(R,4M,J_\Lambda)$ and $(u_n,t_n,c_n)\in \cm(J)$. After extracting a subsequence, 
\begin{itemize}
\itss either there exists a sequence $g\in G=\psl_2(\R)$ such that $(u_n\circ g_n,t_n,c_n)$ converges in $\cm'(J)$, 
\itss or the $u_n$ are vertical discs for $n\gg 1$, $\im u_n=\D\times \{(z'_n,y_n)\}$, with $(z'_n,y_n)\to \partial B_{\infty}(R)$. 
\end{itemize} 
\end{lemma}
  \noindent{\it Proof:} We argue throughout this proof modulo extraction of subsequence. Since $T$ is assumed to be compact, we can assume that $(t_n,c_n)$ converges to $(t,c)\in T\times \cd$ so $J_{(t_n,c_n)}$ converges to $J_{(t,c)}$. Then 
  $(u_n)$ is a sequence of $J_{(t_n,c_n)}$-holomorphic discs with boundaries in $\Gamma_{R,c_n}$, representing the class $E_{c_n}$. By proposition \ref{prop:maincompact}, either the $u_n$ meet $F_{R,M}$ and then are vertical discs so the conclusion of the lemma is achieved, or their images remain in $Y\priv F_{R,M}$ (and obviously in a compact set of $Y$ by maximum principle since $J=J_\st$ at infinity).  In the latter case, $(u_n)$ is a sequence of 
  $J_{(t_n,c_n)}$-holomorphic discs with boundaries in a subset of $\Gamma_{R,c_n}$ contained in a \nbd of a compact subset of $\Gamma_{R,c}$. Their images are moreover contained in a compact subset of $Y$, with bounded  symplectic area. Then Gromov's compactness theorem guarantees the existence of an accumulation point of $(u_n)$ (in the sense of Gromov), which is a non-constant $J_{(t,c)}$-bubble tree. Since however $\pi_2(Y)=0$ and $[u_n]=E_c$ has least area among the non-constant symplectic discs with boundary on $\Gamma_c$, this means that there  is $g_n\in G$ such that $u_n\circ g_n$ converges to an element of $\cm(J)$.\cqfd

We fix  a generic section $J:T\times \cd\to \cj(R,4M,J_\Lambda)$ until the end of this section, and we define 
$$
\begin{array}{l}
W:=\cm(J)\times_G\overline \D:=(\cm(J)\times \overline \D)_{/G},\\
W':=\cm'(J)\times_G\D,
\end{array}
$$
where $g\in G$ acts on $(u,t,c,w)$ by $g\cdot (u,t,c,w):=(u\circ g,t,c,g^{-1}(w))$. We also put 
$$
W^\partial:=W\priv W'=\big(\cm(J)\priv \cm'(J)\big)\times_G\overline \D\, \cup \, \cm'(J)\times_G \partial \D =W^\partial_\text{V}\cup W^{\partial}_\text{H}. 
$$
We now explain that $W$ has the structure of a smooth compact manifold with boundary  $W^\partial$ (and corners). 
Notice first that lemma \ref{le:idiot} obviously shows the following:
\begin{lemma}\label{le:idbdy}
The map \fonction{i}{T\times \{\|(z',y)\|=R\}\times \overline \D}{W^\partial_\text{V}}{(t,z',y,w)}{[t,\re z',w\mapsto(w,z',y,0),w]}
is a one-to one correspondence. 
\end{lemma}
Also, lemma \ref{le:smoothmoduli} provides  $\cm'(J)$ with a structure of  smooth manifold, for which the action of $G$ on $\cm'(J)\times \overline \D$ is smooth, proper and free. As a result, $W'$ is a smooth manifold, $W^\partial_\text{H}$ lies in the closure of $W'$, as a boundary of $W$. 
 
\begin{lemma}\label{le:Wbdy} $W$ is a compact manifold with $W^\partial$ as boundary and corners along $W^\partial_\text H\cap W^\partial_\text V$.
\end{lemma}
\noindent {\it Proof:} 
Let $[t_n,c_n,u_n,w_n]$ be a sequence in $W'$. We can extract a subsequence for which $(t_n,c_n)\to (t,c)\in T\times \cd$. Lemma \ref{le:compactmoduli} leaves us with an alternative, up to extracting a subsequence. Either $\exists g_n\in G$ such that $u_n\circ g_n$ converges to  $u\in \cm'(J)$. We can further extract so that $g_n^{-1}(w_n)$ converges to $w\in \overline \D$, and $[t_n,c_n,u_n,w_n]\to [t,c,u,w]\in W'\cup W^\partial_\text H$. 
 Or the $u_n$ are vertical discs of the form $u_n(z)=(g_n(z),z'_n,y_n')$, where $(z'_n,y'_n)\to \partial B(R)$ and $g_n:\D\to \D$ are ramified coverings. Since the $u_n$ have area $1$ by assumption, these coverings are in fact automorphisms of the discs. Extracting again from $g_n^{-1}(w_n)$ a converging subsequence, we see that $[t_n,c_n,u_n,w_n]$ converge to an element of $W^\partial_\text V$. This discussion also 
provide a parametrization of a \nbd of $W^\partial_\text V$ by $\D\times \{R-\eps<\|(z'_n,y_n')\|\leq R\}$, hence the fact that the points of $W^\partial_\text V\priv W^\partial_\text H$ are boundary points, and the points of $W^\partial_\text V\cap W^\partial_\text H$ are corners.\cqfd

The space $\cm(J)\times_G \overline\D$ comes with a natural evaluation \fonction{\sigma}{W}{T\times Y}{[u,t,c,w]}{(t,u(w)).} As a quotient of a $G$-invariant smooth map on $\cm(J)\times \overline \D$, $\sigma$ is smooth. It  also sends $W^\partial_\text H=\cm(J)\times_G\partial \D$ to $T\times \Gamma_R$ by construction. Notice that lemma \ref{le:compactmoduli} implies that $W^\partial_\text H$ is a smooth manifold with boundary $W^\partial_\text H\cap W^\partial_\text V$. By lemma \ref{le:idbdy}, $\sigma_{|\partial W^\partial_\text H}:\partial W^\partial_\text H\to T\times \partial \Gamma_R$ is a degree $1$ map, hence so is $\sigma_{|W^\partial_\text H}:W^\partial_\text H\to T\times \Gamma_R$. Lemma \ref{le:idbdy} also states that $\sigma_{|W^\partial_\text V}:W^\partial_\text V\to T\times\D\times \{\|(z',y)\|=R,\; x=0\}$ has degree $1$. 
 A consequence of these two points is that $\im \sigma$ is homologous to $T\times Z_R$, relative to $T\times \partial Z_R$. Indeed, the concatenation $(T\times Z_R)\star\im\sigma$, provides an element of $H_{\dim Z_R+\dim T}(T\times Y)$, well-defined because $\partial \sigma$ has degree $1$, and because $T\times Z_R$ and $\im \sigma$ coincide over $T\times \partial Z_R$. But this homology group is 
$$
H_{2k+l+\dim T}(T\times Y)=H_{2k+l+\dim T}(T\times\C^k\times i\R^l\times \T^l)\simeq H_{2k+l+\dim T}(T\times\T^l)=0 \text{ (because $k>0$).}
$$
The image of $\sigma$ is moreover covered by $J_{(t,c)}$-holomorphic discs, for $(t,c)\in T\times \R^k$. Summarizing the discussion of this paragraph, we have obtained:
 \begin{theorem}\label{thm:filling}
 Let $M>0$ be such that the conclusion of proposition \ref{prop:maincompact} holds, and 
 \begin{itemize}
 \itss $Y:=(\C^k\times i\R^l\times \T^l,\om=\pi_*\om_\st)$, \hspace{,5cm}$\cd_R:=\{\|(z',y)\|_\infty\leq R\}\subset \C^{k-1}\times \R^l$,
 
 \itss $\Gamma_R:=S^1\times\cd_R\times \{x=0\}\subset Y$, \hspace{,5cm}  $Z_R=\D\times  \cd_R\times \{x=0\}\subset Y$,
 \itss $F_{R,M}$ an $M$-\nbd of $\D\times \partial \cd_R\times \{x=0\}$ in the $\|(z_1,z',x,y)\|_\infty$-norm,
  \itss $\Lambda\Subset Y\priv (\Gamma_R\cup F_{R,4M})$, \hspace{,2cm} $V_\Lambda\in \op(\Lambda,Y\priv (\Gamma_R\cup F_{R,4M}) )$, \hspace{,2cm} $J_\Lambda:T\times Y \to \cj(\om)$,
 \itss $\cj(R,4M,J_\Lambda):=\{(t,c,J)\in T\times \cd\times \cj(\om)\, |\, J_{|V_\Lambda}\equiv J_\Lambda(t,c),\, J\equiv J_\st \text{ on }\{|(z,y)|>K\}\cup F_{R,4M}\}$ (for some $K\gg 1$).
 \end{itemize}
 For a generic section $J:T\times \cd_R\to \cj(R,4M,J_\Lambda)$ - which exists in any $\cc^\ell$-\nbd of a given smooth  section -, there exists a manifold with boundary and corners $W$ and a smooth map $\sigma:(W,\partial W)\to (T\times Y, T\times \partial Z_R)$ such that:
 \begin{itemize}
 \item[(i)] $[\im \sigma]=[T\times Z(R)]\in H(T\times Y,T\times \partial Z_R,\Z)$,
 \item[(ii)] $\forall p\in \im\sigma$, there exist $(t,c)\in T\times \cd_R$, $u:(\D,\partial \D)\to (Y,\Gamma_{R,c})$ such that 
$p\in \{t\}\times \im u$, $\ca_\om(u)=1$, $\bar \partial_{J(c,t)}u=0$. 
 \end{itemize}
 \end{theorem}
There is nothing to prove here, since all proofs have been done thoughout this section. A remark is however in order. In order to get the point (i), one needs to know the orientability of the moduli spaces, which have not been addressed here. It is a well-known issue in our specific situation, but if one wishes to forget it so as to consider this paper mostly self-contained, the point (i) has to be replaced by 
$$
 [\im \sigma]=[T\times \Z_R]\in H(T\times Y,T\times \partial Z_R,\Z_2).
$$
This is a weaker statement, which is enough for our applications to the $\cc^0$-rigidity of the reduction (see below), and to theorems \ref{thm:loclagcamel} and \ref{thm:loccamel2}, provided the knotting holds in the $\Z_2$-coefficient setting. For theorem \ref{thm:loccamel2} for instance, it amounts to assuming that $\Phi(\R^l\times \{0\})$ intersects $Z_R$ transversally, an odd number of times.

 \subsection{Proof of theorem \ref{thm:loccamel2}}
 Let $\Phi:\R^l\times B^{2n}(a)\to \C^n$ be a parametric symplectic embedding, which satisfies the assumption of theorem \ref{thm:loccamel2} with $C:=4M$, $M$ being the constant obtained in the previous paragraph, in particular theorem \ref{thm:filling}.  Thus, (i) for $|t|\geq K$, $\Phi(t,\cdot)=\id +t$, (ii) $\Phi$ is knotted with $\partial Z_R=S^1\times \cd_R\times \{x=0\}$ (recall that we rescaled the space so that $A=1$), and (iii) $\im \Phi$ does not intersect $F_{R,4M}$. Since $B^{2n}(a)$ is contractible, (ii) amounts to saying that $\im\Phi \cap \partial Z_R=\emptyset$ and $\Phi(\R^l\times \{0\})$ has a non-vanishing homological intersection with $Z_R$, relative to $\partial Z_R$ (say with coefficients in $\Z$, or in $\Z_2$ if one decides to forget about the orientability of the moduli spaces). Restricting $\Phi$ to $\R^l\times B^{2n}(a')$, $a'<a$, we can also assume that $\overline {\im \Phi}\cap \partial Z_R=\emptyset$. 
 
By (i), we can compactify the domain of $\Phi$ to $\T^l\times \overline{B^{2n}(a')}$, getting a parametric symplectic embedding $\Phi:\T^l\times \overline{B^{2n}(a')}\to \C^k\times i\R^l\times \T^l=Y$ which still satisfies (ii), meaning that the homological intersection of $[\Phi(\T^l\times \{0\})]$ with $[Z_R]$ relative to $\partial Z_R$ does not vanish. Since $\dom \Phi$ is now compact, there exists $K>0$ such that $\im \Phi\Subset \{|(z,y)|<K\}$. Taking $T:=\T^l$, we define $\Lambda(t):=\Phi(\{t\}\times \overline{B^{2n}(a')})$ and $J_{\Lambda}(t):=\Phi(t,\cdot)_*J_\st$. Since $\cj(\om)$ is a fiber bundle over $Y$ with contractible fiber, there exists a smooth path $J(t)\in \cj(R,4M,\Lambda(t))$. Applying theorem \ref{thm:filling}, we get a filling $\sigma:(W,\partial W)\to (T\times Y,T\times\partial Z_R)$ with holomorphic discs. Since $[\Phi(T\times \{0\})]\cdot [Z_R]\neq 0$, the graph Graph$(\Phi):=\{(t,\Phi(t,0)),t\in T\}$ has non-vanishing homological intersection with $T\times Z_R$. Moreover, since by theorem \ref{thm:filling} (i), our filling $\im\sigma$ is 
homologous to $T\times Z_R$ relative to $T\times \partial Z_R$, $[\text{Graph}(\Phi)]\cdot[\im \sigma]\neq 0$ and there exists $t\in \T^l$ such that $(t,\Phi(t,0))\in \im \sigma$.   By theorem \ref{thm:filling} (ii), there exists $t'\in \T^l$ and a disc $u:(\D,\partial \D)\to (Y,\Gamma_R)$ such that $(t,\Phi(t,0))\in \{t'\}\times \im(u)$ (hence $t=t'$ and $\Phi(t,0)\in \im(u)$), $\ca_\om(u)=1$, and $\bar \partial_J u=0$ for some $J\in \Gamma^\ell(\cj(M,J_\Lambda))$. This disc is therefore symplectic, has total area $1$, and its trace on the symplectic ball $\Phi(\{t\}\times B^{2n}(a'))$ is $\Phi(t,\cdot)_*J_\st$-holomorphic. The classical monotonicity argument then shows that $1\geq a'$. Since this holds for all $a'<a$, we get $a\leq 1$.\cqfd

  \subsection{Proof of theorem \ref{thm:loclagcamel}}
Let $\Phi:\T^k\times \R^l\hra\C^{k+l}$ be a Lagrangian embedding that satisfies the assumptions of theorem \ref{thm:loclagcamel}: (i) $\Lambda:=\im\Phi$ coincides with $S^1(a)^k\times \R^l$ outside a compact set, and (ii) $L$ is knotted with $\partial Z_R$, meaning that 
$[\Phi(\{*\}\times \R^l)]\cdot[Z_R]\neq 0$, where coefficients belong to  $\Z$ (or $\Z_2$) and homology is understood as always in these paragraphs relative to the boundary. Finally, $\Lambda$ does not meet $F_{R,4M}$.

By (i), we can compactify the domain of $\Phi$ to $\T^k\times \T^l$, getting a Lagrangian embedding $\Phi:\T^n\hra \C^k\times i\R^l\times \T^l=Y
$, which still verifies (ii), and is compactly contained into $\{\|(y,z)\|<K\}$ for some $K\in \R$. We consider $T=\{*\}$, $\Lambda=\im \Phi$, and an almost-complex structure for the moment arbitrary in $\op(\Lambda)$. Applying theorem \ref{thm:filling}, we get for each $J_\Lambda$ a complex structure $J\in \cj(R,4M,J_\Lambda)$ and a $J$-holomorphic disc $u:(\D,\partial \D)\to (Y,\Gamma_c)$ with $[u]=E_c$ and $\im u\cap \Lambda\neq \emptyset$. Moreover, since $J$ is standard at infinity, this disc lies in the compact domain $\{|(y,z)|\leq K\}$, where $K$ does not depend on $J_\Lambda$.

We now proceed to a neck-stretching argument.  
Recall that $\Lambda$ is the image by an exact lagrangian embedding $\Phi$ of $L:=S^1(a)^k\times \T^l$. Consider on $L$ the Euclidean metric $g:=g_\st$. As usual in SFT, we consider a \nbd $\cu$ of $\Lambda$ symplectomorphic to a \nbd of the zero section in $T^*L$, endow it with the metric $g$ induced canonically from $g$ on $L$, and consider $\Sigma:=\{g=\eps_0\}\subset \cu$. For small enough $\eps_0$, this is a contact type hypersurface which splits $Y$ into two pieces $Y^-\cup Y^+$, where $Y^-$ has concave boundary at $X$ (and contains $\partial Z_R$), and $Y^+$ has convex boundary (and contains $\Lambda$). Now we consider a neck-stretching along $\Sigma$. Notice that since $\Lambda$ is disjoint from $\Gamma_R\cup F_{R,4M}$ and contained into $\{|(y,z)|<K\}$, the same holds for $\Sigma$ provided $\eps_0$ is small enough. We can therefore consider a degeneration $J_\eps$ of the complex structure with the following properties. 
\begin{itemize}
\item[\sbull] $J_\Lambda(\eps)$ is defined in $\op(\Lambda)$ and stretches the neck of $\Sigma$ when the parameter $\eps$ goes to $+\infty$. 
\item[\sbull] $J_\eps$ is generic in $\cj(R,4M,J_\Lambda(\eps))$,
\item[\sbull] $J_\eps$ is uniformly bounded in the $\cc^\ell$-topology outside the neck-stretching zone. 
\end{itemize}  
The previous discussion guarantees the existence of a disc $u_\eps:(\D,\partial \D)\to (Y,\Gamma_{c_\eps})$ for some $c_\eps\leq R$, with $\bar \partial_{J_\eps}u_\eps=0$, $[u_\eps]=E_{c_\eps}$,  $u_\eps$ meets $\Lambda$, and $\im u_\eps\subset \{|(y,z)|<K\}$. The compactness theorem in SFT thus implies that there is some $\eps_n$ such that $u_n:=u_{\eps_n}$ converge to a holomorphic building $B$ \cite{boelhowize}, whose main features are summarized below:
\begin{itemize}
\item[\sbull] It has a main component, in $Y^-$, which is a $J_\infty$-holomorphic  map $u_\infty^{\text{main}}:S\to Y^-$, where $S$ is a punctured disc (with $p$ punctures, $p\geqslant 1$ because $u_n$ intersects $\Lambda$, hence $Y^+$), and whose boundary is sent to $\Gamma_R$ and has action $1$. The action is computed with respect to the standard Liouville form $\lambda$ on $Y\simeq\C^k\times T^*\T^l$. 
\item[\sbull] All other components on $Y^-,Y^+$ are $J_\infty$-holomorphic maps whose domains are punctured spheres (because the $u_n$ are discs).
\item[\sbull] The components of the building in $Y^-$ are asymptotic at each puncture to a negative Reeb orbit of $\partial Y^-=\Sigma$ (the boundary of the image is oriented by the opposite of the Reeb flow). Similarly, the components in $Y^+$ are asymptotic to positive Reeb orbits of $\Sigma$. 
\item[\sbull] There might be intermediate layers: components of $B$ in $\Sigma\times \R$ (the symplectization of $\Sigma$). These components are $J$-holomorphic for some cylindrical almost complex structure, again punctured spheres, asymptotic to positive Reeb orbits at $\Sigma\times\{+\infty\}$ and negative Reeb orbits at $\Sigma\times \{-\infty\}$. 
\item[\sbull] The total symplectic area of these components is $1$, and they glue together to form a topological disc (see figure \ref{fig:sft}). We will refer to this last property by saying that the building $B$ forms a disc.  
 \end{itemize} 
 
 \begin{figure}[h!]
\begin{center}
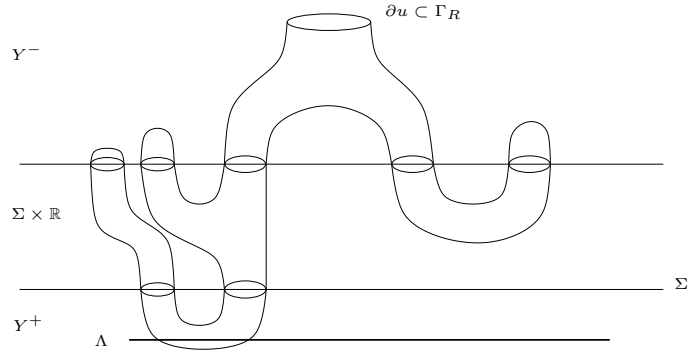
\end{center}
\caption{A holomorphic building.}\label{fig:sft}
\end{figure}

 Denote by $-\gamma_1,\dots,-\gamma_p$ the (negative) Reeb orbits to which $u_\infty^{\text{main}}$ is asymptotic. Since  $B$ forms a disc, its other components glue together to form $p$ discs asymptotic to $\gamma_1,\dots, \gamma_p$. Since the $\gamma_i$ are Reeb orbits, they project to closed geodesics  of $L$ under the natural projection $\pi:Y^+\simeq T^*_{\eps_0}L\to L$ (recall that $Y^+\simeq \{g<\eps_0\}$). We write $[\pi\gamma_i]=\sum k_j^{(i)}e_j\in H_1(L)$ for the homology classes of these geodesics, where $e_j$ is the class of the $j$-th $S^1$-factor in $L\simeq (S^1)^n$, and $k_j^{(i)}\in \Z$. Observe at this point that since $\gamma_i$ bounds a topological disc, $k_j^{(i)}=0$ for $j> k$ (because $Y\simeq  \C^{k}\times T^*\T^l$).  We will use the following obvious fact:
\begin{fact} If $\gamma$ is a positive Reeb orbit of $\Sigma$, the cylinder \fonction{\rho_\gamma}{[0,1]\times S^1}{Y^+\simeq T^*_{\eps_0}L}{(s,t)}{s\gamma(t)} is symplectic, with (oriented) boundaries $\gamma$ and $-\pi\gamma$. 
\end{fact}

\begin{figure}[h!]
\begin{center}
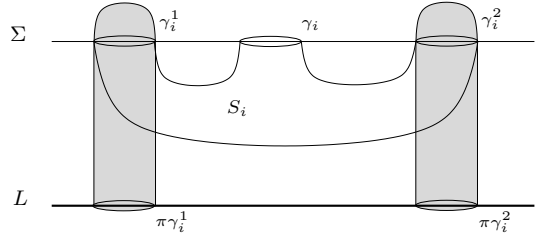
\end{center}
\caption{Capping the class $[\pi\gamma_i]$.}\label{fig:cap}
\end{figure}

\begin{lemma}
For each $i$, $\sum k_j^{(i)}>0$. 
\end{lemma}
\noindent {\it Proof:}  Denote by $S_i$ the connected component of $B$ in $Y^+$ which is glued to $u_\infty^{\text{main}}$ along $\gamma_i$. This component cannot be a disc because $\pi\gamma_i$ is non-contractible in $Y^+\subset T^*\T^n$. Therefore, $S_i$ has boundary components asymptotic to positive Reeb orbits $\gamma_i,\gamma_i^1,\dots,\gamma_i^{k_i}$, with $k_i\geqslant 1$. Observe that $\pi S_i$ provides a chain whose boundary is $[\pi\gamma_i]+\sum [\pi\gamma_i^k]$, so $[\pi\gamma_i]=-\sum [\pi\gamma_i^k]$.  But since $B$ forms a disc, there is a disc with positive area - composed of a gluing of different components of $B$ -
which is asymptotic to each $-\gamma_i^k$. Gluing to these discs the cylinders $\rho_{\gamma_i^k}$, we get discs with positive area,  whose boundaries lie in $L$ and represent the class $-[\pi\gamma_i^k]$ (see figure \ref{fig:cap}). Therefore, 
$$
\hspace{4cm}0<\sum \lambda[-\pi\gamma_i^k]=\lambda[\pi\gamma_i]=\sum_{i=1}^{k} k_j^{(i)} \, a.\hspace{4,5cm} \square
$$

The map obtained by gluing to $u_\infty^\text{main}$  the symplectic cylinders $\rho_{\gamma_i}$ now gives a symplectic surface of area 
$$
0<1-\sum \lambda[\pi\gamma_i]=1-\sum_{i=1}^p\left(\sum_{j} k_j^{(i)}\right)a\leqslant 1-a,
$$
where the last inequality holds because $p\geq 1$ and the previous lemma. We get $1\geq a$.\cqfd

\section{Coisotropic $\cc^0$-rigidity}\label{sec:coisorig}
\subsection{The reduction of a coisotropic submanifold}\label{sec:defred}

We now define precisely what we mean by reduction of a coisotropic submanifold. This is just a generalization of the discussion for hypersurfaces, which was held in \cite{buop}. 

\begin{definition} \label{D:reduction}

Let $ \Sigma^m $ be a smooth manifold endowed with a foliation $ \mathcal{F} $, where the dimension of each leaf of the foliation equals $ r $.
\begin{enumerate}

\item
We say that an open topological submanifold$\,$\footnote{Recall that a topological submanifold of a topological manifold $ X $ is a subset $ Y \subset X $, such that there exists a topological manifold $ Z $ and a map $ i : Z \rightarrow X $ which is a homeomorphism onto the image $ i(Z) = Y $.} $ \; U^{m-r}\subset \Sigma$ is (topologically) transverse to the foliation $ \mathcal{F} $ on $ \Sigma $ if $ U $ has a \nbd $ V\subset \Sigma$ such that $ U $ intersects exactly once each leaf of $ V $ (where by a leaf of $ V $ we mean a connected component of the intersection of $ V $ with a leaf of $ \mathcal{F} $). 

\item 
Let $ U, U' \subset \Sigma$ be $ (m-r) $-dimensional topological submanifolds, that are transverse to the  foliation $ \cf $. We say that $ U $ and $ U' $ are equivalent (denoting $ U\sim U'$) if there exists a (continuous) homotopy $ G : W \times [0,1] \rightarrow \Sigma $, $ t \in [0,1] $, of a topological manifold $ W^{m-r} $, such that $ G_{|W \times \{ 0 \}} $ is a homeomorphism onto $ U $, $ G_{|W \times \{ 1 \}} $ is a homeomorphism onto $ U' $, and such that for each $ x \in W $, the trajectory $ t \mapsto G(x,t) $ goes along a leaf of $ \mathcal{F} $.

\item
The reduction of a smooth manifold $\Sigma$ endowed with a foliation $ \mathcal{F} $, denoted by $\red(\Sigma,\mathcal{F})$, is defined as the set of open topological submanifolds $U^{m-r} \subset \Sigma $ which are transverse to the characteristic foliation of $ \Sigma $, considered modulo the above equivalence relation. If $ \Sigma $ is a coisotropic (or, more generally, a pre-symplectic) submanifold of a symplectic manifold $ M $, and $ \mathcal{F} $ is the characteristic foliation on $ \Sigma $, then we simplify the notation for the reduction to $ \red(\Sigma) $.

\end{enumerate}
\end{definition}

\medskip

\noindent Now let us address several points:

\paragraph{Smooth and symplectic structures} \label{D:sss} On a topological submanifold $ U \subset \Sigma $ which is transverse to the foliation $ \mathcal{F} $, we have a natural structure of a smooth manifold. Moreover, if $ \Sigma $ is a coisotropic (or, more generally, a pre-symplectic) submanifold of a symplectic manifold, and $ \mathcal{F} $ is the characteristic foliation on $ \Sigma $ (we will call such situation a ``symplectic setting"), then $ U $ moreover inherits a natural symplectic structure. Indeed, let $ V $ be a neighbourhood of $ U $ in $ \Sigma $ such that $ U$ intersects exactly once each characteristics of $ V $, as in definition~\ref{D:reduction}. Then any point $ z \in U $ has a neighbourhood $ U_1 \subset U $ such that $ U_1 $ lies inside a (smooth) flow-box $ \Phi : W_1 \times (0,1)^r \rightarrow \Sigma $, where $ \im \Phi \subset V $. Then the map $ \phi := \pi \circ \Phi^{-1} : U_1 \rightarrow W_1 $ is injective and hence, by the Invariance of Domain theorem, is a homeomorphism onto the open image $ \phi(U_1) \subset W_1 $ (here $ \pi : W_1 \times (0,1)^r \rightarrow W_1 $ is the natural projection). The map $ \phi $ induces a natural smooth structure on $ U_1 $, in case of a general foliation, and moreover induces natural symplectic structure in $ U_1 $ if we are in a symplectic setting.

\paragraph{Naturality of the structures} If two topological submanifolds $ U, U' \subset \Sigma$ are equivalent ($ U\sim U'$), then they are diffeomorphic (and moreover symplectomorphic if we are in a symplectic setting) via a homotopy $ G : W \times [0,1] \rightarrow \Sigma $, as in definition \ref{D:reduction}. Let us describe explicitly the diffeomorphism (resp. symplectomorphism) between $ U $ and $ U' $. By continuity of $ G $ and since $ U $ is topologically transverse to leaves of the foliation $ \mathcal{F} $, for any point $ z \in W $ and any $ t \in [0,1] $ there exists a neighbourhood $ W_1 \Subset W $ of $ z $ such that the closure of the image $ G(W_1 \times \{ t \}) $ lies inside a (smooth) flow-box $ \Phi : W_2 \times (0,1)^r \rightarrow \Sigma $, and moreover the map $ \phi_t := \pi \circ \Phi^{-1} \circ G : W_1 \times \{ t \}  \rightarrow W_2 $ is injective (here $ \pi : W_2 \times (0,1)^r \rightarrow W_2 $ is the natural projection, as before). Then, by the Invariance of Domain theorem, $ \phi_t $ is a homeomorphism onto the open image $ W_3 := \phi_t (W_1) \subset W_2 $. This induces a smooth (resp. symplectic) structure on $ W_1 \times \{ t \} $. Moreover, since $ W_1 \Subset W $ and $ G(W_1 \times \{ t \}) \Subset \Phi (W_2 \times (0,1)^r) $, it follows that we also have $ G(W_1 \times \{ t' \}) \Subset \Phi (W_2 \times (0,1)^r) $ as well and moreover $ \phi_{t'}(z,t') = \phi_t(z,t) $ for every $ z \in W_1 $, whenever $ t' \in [0,1] $ is sufficiently close to $ t $ (here $ \phi_{t'} = \pi \circ \Phi^{-1} \circ G : W_1 \times \{ t' \}  \rightarrow W_2 $). Hence the induced smooth (resp. symplectic) structures on $ W_1 \times \{ t \} $ and on $ W_1 \times \{ t' \} $ coincide, when $ t' $ is sufficiently close to $ t $.

\paragraph{The induced map} Let $ h: \Sigma \rightarrow \Sigma' $ be a homeomorphism between smooth manifolds $ \Sigma $ and $ \Sigma' $ which are endowed with foliations $ \mathcal{F} $ and $ \mathcal{F}' $, such that $ h $ maps $ \mathcal{F} $ to $ \mathcal{F}' $. Then $ h $ defines a natural map $ \hat{h} : \red(\Sigma,\mathcal{F}) \rightarrow \red(\Sigma',\mathcal{F}') $ by $ \hat{h}([U]) := [h(U)] \subset \Sigma' $. Clearly, the definition does not depend on the representative $ U $ of $ [U] $. Recall that by a theorem of Humili\`ere, Leclercq and Seyfaddini~\cite{hulese}, if $ h : M \rightarrow M' $ is a symplectic homeomorphism and if $ h $ maps a smooth coisotropic submanifold $ \Sigma $ onto a smooth submanifold $ \Sigma' $, then $ \Sigma' $ is coisotropic, and the restriction $ h_{|\Sigma} : \Sigma \rightarrow \Sigma' $ preserves the characteristic foliation, and hence we get a natural induced map $ \hat{h} : \red(\Sigma) \rightarrow \red(\Sigma') $. \\

\subsection{Proof of theorem \ref{thm:rigredloc}}
Let us adopt the following notation for the sake of fluidity. Given two submanifolds $A,B$ with boundaries, we say that $A$ is knotted with $\partial B$, relative to $\partial A$, if any homotopy of $A$ relative to $\partial A$ that avoids $\partial B$ intersects $B$. Notice that here, the filling $B$ of $\partial B$ is fixed, and not made explicit in the sentence. In the following still, we will be interested in a knotting with $\partial Z(A,R)$, whose filling is implicitely meant by $Z(A,R)$.

Let $h$ be a symplectic homeomorphism defined in a \nbd $U$ of $B^{2k}(1)\times (-1,1)^l$ in $\C^{k+l}=\C^n$, with values in $\C^n$, such that  $h\big(B^{2k}(1)\times (-1,1)^l\big) \subset \C^k\times \R^l$. We also assume for later convenience that $U$ is contractible. We need to prove that for $\delta$ small enough, for any $a<\delta$, if $h(B^{2k}(a)\times (-1,1)^l)\subset Z(A)\times \R^l$, then $A\geq a$. 

Since $h$ is a symplectic homeomorphism, we know by \cite{hulese} that it takes the characteristic leaves $\{*\}\times (-1,1)^l$ into $\{*\}\times \R^l$. In other terms, there exists a continuous function $\hat h:B^{2k}(1)\to \C^k$ and  an open subset $\Om(z)\subset \R^l$ for each point $z\in B^{2k}(1)$ such that 
$$
h(\{z\}\times (-1,1)^l)=\{\hat h(z)\}\times \Om(z). 
$$
Since $Z(A)$ is invariant by translation along the $z'$-axis, we can assume that $\hat h(0)=(z_1,0)$, which we do henceforth (recall that we split $\C^k=\C\times \C^{k-1}$, with  $(z_1,z')$ the split coordinates). We divide this proof in four steps:

\noindent\underline{Step 1: Adjusting $\hat h$ and $\Om$.} As already observed in section \ref{sec:defred},  the map $\hat h$ is a local homeomorphism. Thus, its restriction to $B^{2k}(R)$ for $R$ small enough is injective, so after maybe restricting $h$ to $B^{2k}(R)\times (-1,1)^l$, we can assume that $\hat h$ is injective, which we do henceforth. 
Moreover, since $h$ is a homeomorphism, the map $z\mapsto \Om(z)$ is continuous, so by a further restriction, we can also ensure that 
$$
\Om(z)\supset \frac 12\Om(0) \hspace{1cm} \forall z \in B^{2k}(R). 
$$
Finally, since $\Om(0)$ is homeomorphic to $(-1,1)^l$, there exists a compactly supported diffeomorphism $\phi:\R^l \to \R^l$  such that $\phi(\Om(0))$ is arbitrarily close to $(-4,4)^l$. Letting $\Phi:T^*\R^l\to T^*\R^l$ the natural lift of $\phi$ and considering $(\id \otimes \Phi)\circ h$ instead of $h$, we can therefore also assume that $\Om(0)$ is close to $(-4,4)^l$, and that 
$\Om(z)\supset (-2,2)^l$ for all $z$. Summarizing this discussion, we see that we can assume without loss of generality that 
\begin{itemize}
\its $h$ is a symplectic homeomorphism defined in a \nbd of $B^{2k}(1)\times (-1,1)^l$ with values in $\C^n$, that takes $B^{2k}(1)\times (-1,1)^l$ to a subset 
$$\{(z,t)\in W\times \R^l,\; t\in \Om(z)\},$$
\its $\hat h:B^{2k}(1)\to \C^k$ is injective, 
\its $ \Om(z)\supset (-2,2)^l$ for all $z$. 
\end{itemize}
Obviously, $W$ is open, and if we fix any compact exhaustion $W_n$ of $W$, there exists therefore $\delta_n>0$ such that 
$$
\im h\supset Q_n:=W_n\times (-1,1)^l\times (-\delta_n,\delta_n)^l. 
$$
Moreover, since $h$ is a homeomorphism, if $\delta_n$ is chosen small enough, 
$$
\im h\cap Q_n\cap \C^k\times \R^l=h(B^{2k}(1)\times (-1,1)^l). 
$$
In other terms, $h$ intersects the coisotropic subset $\C^k\times \R^l$ in $Q_n$ exactly along $B^{2k}(1)\times (-1,1)^l$. 

\noindent\underline{Step 2: Knotting $\hat h$.} Let $n_0$ be such that $\hat h(0)\notin W_{n_0}$, and put $R_0:=\frac 12d(\hat h(0),\partial W_{n_0})$ (measured in the infinity norm). For $\delta_1(h)\ll 1$, $d\big(\hat h(\partial B(\delta_1(h))),\hat h(0)\big)<R_0$. Since $\hat h(0)$ has the form $(z_1,0)$, if $a\leq \delta_1(h)$ and $\hat h(B(a))\subset Z(A)$, $\hat h(B(a))\subset D(A)\times B_\infty(R_0)$. Moreover, the first step above guarantees that 
$$
h(\{0\}\times (-1,1)^l)\cap D(A)\times B_\infty(R_0)\times B_\infty(\delta_n)\cap Q_{n_0}=\{\hat h(0)\}\times (-1,1)^l
$$
Thus, provided that $R<R_0$ and $Z(A,R)\subset Q_{n_0}$, $h(B(a)\times (-1,1)^l)$ is knotted with $\partial Z(A,R)$ in the sense given at the beginning of this section. This condition can be written:
\begin{equation}\label{eq:cond1}
R<R_0,\hspace{,3cm} A\leq R_0,\hspace{,3cm}  R<\delta_{n_0}. 
\end{equation}

Notice also that provided $A$ is chosen small enough compared with $R$, $h\big(B(a)\times (-1,1)^l\big)$ is even knotted with $\partial Z(A,R)\cup F(A,R,CA)$, where the filling of this last set is still by $Z(A,R)$. For this to hold, it is necessary and sufficient that $F(A,R,CA)\cap \C^k\times \R^l=\emptyset$, which holds as soon as 
\begin{equation}\label{eq:cond2}
R>CA,
\end{equation}
and $F(A,R,CA)\subset Q_{n_0}$, which needs in top of \eqref{eq:cond1}
\begin{equation}\label{eq:cond3}
CA<\delta_{n_0}. 
\end{equation}
To conclude this step, notice that the inequalities \eqref{eq:cond1}-\eqref{eq:cond3} can be achieved as soon as 
$$
A<\delta_2(h):=\min(1,\frac 1C)\min(R_0,\delta_{n_0})
$$
because then $R$ can be chosen between $CA$ and $\min(R_0,\delta_{n_0})$ and the inequalities are then all satisfied. 

\noindent\underline{Step 3: from $h$ to a symplectomorphism.}   Let now assume that $A<\delta(h):=\min(\delta_1(h),\delta_2(h))$. Since $h$ is a symplectic homeomorphism, it can be approached in the uniform norm by symplectic diffeomorphisms $f_n:U\to \C^n$. Using the contractibility assumption on $U$, it is easy to see that $f_n$ can be extended to a compactly supported Hamiltonian diffeomorphism of $\C^n$. Moreover, since $f_n$ is close to $h$ in $U$, the following two points are satisfied for large enough $n$: 
\begin{itemize}
\its $f_n\big(B(a)\times (-1,1)^l\big)$ is knotted with $Z(A,R)\cup F(A,R,CA)$ in $Q_{n_0}$,
\its $f_n(U)\supset Q_{n_0}$, so $f_n\big(B(a)\times \R^l\priv(-1,1)^l\big)\cap Q_{n_0}=\emptyset$. 
\end{itemize} 
As a result, and since $Z(A,R)\cup F(A,R,CA)\subset Q_{n_0}$, $f_n\big(B(a)\times \R^l\big)$ is knotted with $Z(A,R)\cup F(A,R,CA)$ so corollary \ref{cor:camelcoiso} guarantees that $A\geq a$. 

\noindent\underline{Step 4: from $A$ to $a$.} We have therefore got our non-squeezing inequality when $A<\delta(h)$. It only remains to get it for $a<\delta(h)$ in order to establish our statement. But this is obvious. Indeed, 

if $a<\delta(h)$, either $A\geq \delta(h)$, and then $A>a$, or $A<\delta(h)$ and then the previous analysis implies that $a\geq A$. \cqfd

\subsection{Proof of theorem \ref{thm:rigred2}}
Theorem \ref{thm:rigred2} can be seen as a corollary of theorem \ref{thm:rigredloc}. Indeed, what we have to prove is that if $\Sigma,\Sigma'$ are two $n+1$-dimensional coisotropic submanifolds  of $2n$-dimensional symplectic manifolds, if there exists a symplectic homeomorphism $h:(\op(\Sigma),\Sigma)\lra (\op(\Sigma'),\Sigma')$, the reduction $\hat h$ of $h$ is area preserving. Since area preservation is a local property, it is enough to check it for small elements of the reduction. We can therefore assume that both $\Sigma$ and $\Sigma'$ are $D^2(1)\times \R^{n-1}\subset \C^n$. We now argue by contradiction. Assume that $\hat h$ is not area preserving, and considering maybe $h^{-1}$ instead of $h$, we can assume that it expands the area of some elements of the reduction. Then there exists a subdisc $D(R_0)\Subset D(1)$ whose area is increased by  $\hat h$. Moreover, there must exist a point $z\in D(R_0)$ and a decreasing sequence $\eps_n>0$ that tends to $0$ such that the area of $\hat h(D(z,\eps_n))$ exceeds $\eps_n$ for all $n$. But this is precisely ruled out by theorem \ref{thm:rigredloc}. \cqfd

{\footnotesize
\bibliographystyle{alpha}
\bibliography{biblio}
}

\paragraph{Emmanuel Opshtein,} $ $\\
Institut de Recherche Math\'ematique Avanc\'ee\\
UMR7501, Université de Strasbourg et CNRS\\
7 rue Ren\'e Descartes\\
67000 Strasbourg, FRANCE.\\
E-mail: opshtein@unistra.fr 

\end{document}